\documentclass[12pt,final,a4paper,titlepage]{article}
\usepackage{amssymb}
\usepackage{layout}

\newtheorem{theorem}{Theorem}[section]
\newtheorem{corollary}{Corollary}[section]
\newtheorem{conjecture}{Conjecture}[section]
\newtheorem{lemma}{Lemma}[section]
\newtheorem{definition}{Definition}[section]
\newtheorem{remark}{Remark}[section]

\makeatletter
\@addtoreset{equation}{section}
\renewcommand{\p@equation}{}
\makeatother

\begin{document}

\begin{titlepage}


\begin{center}
{\Large\textbf{Regularity of dissipative differential operators}}
\end{center}

\bigskip

\bigskip

\begin{center}
\textbf{A.Minkin}
\footnotetext{
Partially supported by the International Soros Science Education Program,\\
grant N 149d "Soros Associate Professor"}

\end{center}

\bigskip

\bigskip

\hspace{4cm} Novouzenskaya 51/63, Appt. 84,

\medskip

\hspace{4cm} Saratov, 410017, Russia

\end{titlepage}

\thispagestyle{empty}

\bigskip

\noindent \textbf{Mailing address:} Prof. A.Minkin

\medskip

\hspace{4cm} Novouzenskaya 51/63, Appt. 84,

\medskip

\hspace{4cm} Saratov, 410017, Russia

\medskip

\hspace{4cm} e-mail: mariya@rakhmat.saratov.su
\newpage

\thispagestyle{empty}

\begin{abstract}{\footnotesize
S.G.Krein's conjecture concerning Birkhoff-regularity of dissipative differential
operators has been proved in the even order case. As a byproduct an existence of
the limit of characteristic matrix as $\lambda \rightarrow \infty $ in the
lower half-plane has been established. Up to multiplication by a nonvanishing
matrix ${D}$
this limit coincides with the ratio
\[
\mathbf{\Theta}(b^0,b^1)^{-1}
\cdot
\mathbf{\Theta}(b^1,b^0)\cdot {D}
\]
of the matrices of regularity determinants.

}
\end{abstract}

\hfill\emph{To A.P.Khromov with admiration and affection}

\vspace{1cm}

\section{Introduction}

\subsection{History}

Spectral theory of boundary value problems in a finite interval is one of
the most elaborated parts of the theory of linear nonselfadjoint operators.
Just at the beginning of the century G.D.Birkhoff discovered a class of the
so-called \textit{regular }boundary conditions (further on,
Birkhoff-regular) with a lot of remarkable spectral properties: estimate of
the Green's function, asymptotics of eigenvalues and eigenfunctions, point
convergence similar to that of trigonometric Fourier series and so on \cite
{Nai:ldo}.

However, further investigations revealed absence of \textit{concrete classes}
of boundary conditions with the same list of properties if
Birkhoff-regularity is violated. In this case the resolvent admits a
polynomial and even an exponential growth. More precisely, for separated
boundary conditions the Green's function grows exponentially at least in one
of the triangles $0<x<t<1$ or $0<t<x<1$,  and in 1991 A.P.Khromov described
a class of third order irregular problems \cite{hro91a} where it has an
exponential growth in both triangles.

From the other hand the general theory of linear operators provides such
\textit{good} subsets as self-adjoint, normal or dissipative ones. Of
course, the list may be increased but we stop here. Quite naturally there
raises a question of an intersection between them and the concrete class of
differential operators generated by Birkhoff-regular boundary value
problems. Recall first that \textit{regularity} of self-adjoint boundary
conditions has been proved by S.Salaff in 1968 for even order operators \cite
{sal68}, see also work of H.Fiedler \cite{fie72}. In the odd case this fact
was established in 1977 \cite{min77a}. Later these results were generalized
in \cite{min93a}.

At the beginning of nineties we supposed validity of the following
conjecture.

\begin{conjecture}
Dissipative differential operators are Birkhoff-regular.
\end{conjecture}

\noindent Recently V.A.Il'in has kindly informed us that it belongs to
S.G.Krein, to one of his talks during the Voronez mathematical schools in
seventies-eighties (personal communication). Our main result is as follows.

\begin{theorem}
\label{thm0.1}Even order dissipative differential operators are
Birkhoff-regular.
\end{theorem}

In the odd order case our demonstration provides nonvanishing of only one
regularity determinant instead of the two ones involved in the definition (%
\ref{eq1.5}) of Birkhoff-regularity.

\subsection{Abstract approach}

Before passing to the proof it should be pointed out that the spectral
theory of abstract dissipative operators is deeply explored by means of the
functional model theory \cite{SNF:harm,Nik:shift}. For instance, S.R.Treil'
has found a strong criterion of unconditional basicity of eigenfunctions
\cite[theorem 14.1]{tre:89}. Observe that it requires their uniform
minimality, which seems to be unattainable unless $L\in (R)$. Hence, the
abstract approach turns out to be ineffective for {boundary value problem}s.%
\hfill

Instead, an \textit{almost orthogonality} property of the Birkhoff's
fundamental system of solutions serves below as a main tool. It has been
proved in \cite{min80} for ordinary differential expressions, in \cite
{min93b} for the quasidifferential ones with a summable coefficient by
the $(n-1)$-st derivative and asserts that\hfill
\begin{equation}  \label{eq0.1}
\left\| \sum_{k=0}^{n-1}c_k\cdot y_k(x,\varrho ) \right\| _{L^2(0,1)}^2
\asymp \sum_{k=0}^{n-1} \left| c_k\right| ^2\cdot \left\| y_k(x,\varrho
)\right\| _{L^2(0,1)}^2
\end{equation}
for any constants $c_k$ which may vary with $\varrho $. Other applications
of (\ref{eq0.1}) may be found in \cite{mishu97}.

\subsection{Notations}

Throughout, matrices with block entries are written in boldface, for
instance, $\mathbf{\Delta }=[\Delta_{jk}]$ stands for the matrix of a
determinant $\Delta $. Components of matrices and vectors are enumerated
beginning from zero, for instance, '$0$-th' row, '$k$-th' column etc. We use
few abbreviations:

\begin{itemize}

\item  $a:=b$ or $d=:c$ means that $a$ equals $b$ or $c$ equals $d$ by
       definition;

\item  $[a]:=a+O(1/\varrho )$ stands for the Birkhoff's symbol;

\item  $A\asymp B$ means a double-sided estimate $\ \ C_1\cdot |A|\leq
|B|\leq C_2\cdot |A|\ \ $ with some absolute constants $C_{1,2}$, which
don't depend on the variables $A$ and $B$. In this case we shall say that $A$
\emph{is equivalent to} $B$;

\item  $\mathbb{C}_{\pm }$ - upper/lower half-plane, $\mathbb{R}$ stands for
the real axis.
\end{itemize}

Different constants are denoted $C, C_1, c$ and so on. They may vary even
during a single computation.

\subsection{General remarks}

The paper is organized as follows. In \S\ \ref{sec:prelim} we recall the
definition of Birkhoff-regularity and build a suitable form of the Green's
function. In \S\ \ref{sec:bound} boundedness of a characteristic matrix is
established. The proof is completed in \S\ \ref{sec:end}. Then at the end of
the paper we give some remarks concerning applications and possible
generalizations.

\section{Preliminaries\label{sec:prelim}}

\subsection{Birkhoff-regular problems}

Consider a differential operator $L$ in $L^2(0,1)$ defined by a two-point
boundary value problem ($(D=-id/dx)$):\hfill
\begin{equation}
l(y)\equiv D^ny+\sum_{k=0}^{n-2}p_k(x)D^ky=\lambda y,\ \ 0\le x\le 1,\ \
p_k\in L(0,1)  \label{eq1.1}
\end{equation}
and $n$ linearly independent normalized boundary conditions \cite[p.65--66]
{Nai:ldo}:
\begin{equation}
U_{\jmath }(y)\equiv b_{\jmath }^0D^jy(0)+b_{\jmath }^1D^jy^{(}1)+\ldots
=0,\;\;j=0,\ldots ,n-1.  \label{eq1.2}
\end{equation}
\noindent
Here the ellipsis takes place of lower order terms at $0$ and at $1$.
Further $b_{\jmath }^0,b_{\jmath }^1$ are column vectors of length $%
r_{\jmath }$, where
\[
0\le r_{\jmath }\le 2,\;\sum_{k=0}^{n-1}r_{\jmath }=n,\ \
rank\left(b_{\jmath }^0b_{\jmath }^1\right)=r_{\jmath }.
\]
Such form of normalized boundary conditions was introduced by S.Salaff
\cite[p.356--357]{sal68}. It is evident that $r_{\jmath }=0$ implies the
absence of order $j$ conditions. In the case $r_{\jmath }=2$ we merely put
\[
\left( b_{\jmath }^0b_{\jmath }^1\right) =\left(
\begin{array}{cc}
1 & 0 \\
0 & 1
\end{array}
\right) .
\]


\begin{definition}
Let $q=Entier(n/2)$, $\varepsilon _{\jmath }:=\exp (2\pi ij/n)$, $k=0,..,n-1;
$
\begin{equation}
b^i=(b_{\jmath }^i)_{j=0}^{n-1}=\left(
\begin{array}{l}
{b_0^i} \\
\vdots  \\
b_{n-1}^i
\end{array}
\right) ,\;B_k^i=\left( b_{\jmath }^i\cdot \varepsilon _k^j\right)
_{j=0}^{n-1},\;i=0,1;  \label{eq1.3}
\end{equation}
\medskip
\begin{equation}
\Theta _p(b^0,b^1)=\left| B_k^0,k=0,\ldots ,p-1|B_k^1,k=p,\ldots ,n-1\right|
,  \nonumber  \label{eq1.4a}
\end{equation}
\begin{equation}
\Theta (b^0,b^1):=\Theta _q(b^0,b^1)  \label{eq1.4}
\end{equation}
The vertical line $|$ separates columns with superscripts $0$ and $1$. The
determinant $\Theta _p(b^1,b^0)$ is defined like (\ref{eq1.4}) by
interchanging the superscripts $0$ and $1$. We shall call boundary
conditions (\ref{eq1.2}) and the corresponding operator $L$ Birkhoff-regular
and write $L\in (R)$ if\hfill
\begin{equation}
\left\{
\begin{array}{llll}
\Theta (b^0,b^1) & \neq  & 0, & n=2q, \\
\Theta (b^0,b^1) & \neq  & 0\;\mbox{and }\;\Theta (b^1,b^0)\neq 0, & n=2q+1.
\end{array}
\right.   \label{eq1.5}
\end{equation}
\end{definition}

\begin{definition}
Birkhoff-regular {boundary value problem} is \emph{strongly regular} if
either $n$ is odd or if it is even, $n=2q$, and the polynomial\hfill
\[
\ F(s):=\left| B_0^0+s\cdot B_0^1,\,B_k^0,k=1,\ldots ,q-1\,|\,B_q^1+s\cdot
B_q^0,\,B_k^1,k=q+1,\ldots ,n-1\right|
\]
has two simple roots.
\end{definition}

\noindent
This form of Birkhoff-regularity was invented by S.Salaff \cite[p.361]{sal68}
who has done a first serious investigation of the nature of the regularity
determinants. The reader must keep in mind that our determinants differ
slightly from those in \cite[p.361]{sal68}, namely $\Theta (b^1,b^0)$
coincides with $\Theta (b^0,b^1)$ from \cite{sal68} but this does not affect
further considerations.


\subsection{Particular solution}

Set
\(
\varrho =\lambda ^{1/n},\quad |\varrho |=|\lambda |^{1/n}
\)
and\hfill
\begin{equation}
\arg \varrho =\arg \lambda /n,\qquad 0\leq \arg \lambda < 2\pi .
\label{eq1.6}
\end{equation}
Then $\varrho \in S_0\cup S_1$, where
\[
S_k=\{\varrho \,\,\,|\,\,\,\pi k/n\le \arg \varrho < \pi (k+1)/n\}.
\]
Note that in every sector $S_k$ there exists a fundamental system of
solutions $\{y_{\jmath}(x,\varrho )\}_{j=0}^{n-1}$ with an exponential
asymptotics:\hfill
\begin{equation}
D^ky_{\jmath }(x,\varrho )=(\varrho\varepsilon _{\jmath })^k\cdot \exp
(i\varrho \varepsilon _{\jmath }x)[1],\ \ j,k=0,\ldots ,n-1.  \label{eq1.7}
\end{equation}
In the sequel it will be convenient to introduce a number $p$ such that
solutions $y_{\jmath }(x,\varrho )$ decay as $j<p$ and exponentially grow
otherwise. Clearly, $p$ depends upon the sector's choice and its values are
presented in the table \ref{table:p}.

\begin{table}[h]
\hfil
\begin{tabular}{|l|l|l|}
\hline
$\varrho \diagdown n$ & $2q$ & 2q+1 \\ \hline
$\in S_0$ & $q-1$ & $q$ \\
$\in S_1$ & $q-1$ & $q-1$ \\ \hline
\end{tabular}
\hfil
\caption{Values of $p$}
\label{table:p}
\end{table}
Introduce a wronskinian\hfill
\[
W(x,\varrho )=\left| D^ky_{\jmath }\right| _{j,k=0}^{n-1}
\]
and let $W_{\jmath }(x,\varrho )$ be the algebraic complement of the element
$D^{n-1}y_{\jmath }$. Set
\(
\tilde y_{\jmath }(x,\varrho ):=W_{\jmath }/W.
\)
It is easy to check that\hfill
\begin{equation}
\tilde y_{\jmath }(x,\varrho )=\frac 1{n(\varrho \varepsilon _{\jmath
})^{n-1}}\exp (-i\varrho \varepsilon _{\jmath }x)[1].  \label{eq1.8}
\end{equation}
Introducing the kernel\hfill
\[
g_0(x,\xi ,\varrho )=i\cdot \left\{
\begin{array}{rl}
\sum\limits_{k=0}^{p-1}\varepsilon _k^{-(n-1)}y_k(x,\varrho) \tilde
{y_k}(\xi,\varrho ),\ \ x>\xi &  \\
-\sum\limits_{k=p}^{n-1}\varepsilon _k^{-(n-1)}y_k(x,\varrho) \tilde
{y_k}(\xi,\varrho ),\ \ x>\xi &
\end{array}
\right.
\]
we get a particular solution $g_0(f)$ of the equation $l(y)=\lambda y+f$,
\begin{equation}
g_0(f):=\int\limits_0^1g_0(x,\xi ,\varrho )f(\xi ))d\xi .  \label{eq1.9}
\end{equation}

\subsection{New fundamental system of solutions}

In the sequel it will be more convenient to use another fundamental system
of solutions $\{z_k\}_{k=0}^{n-1}$, where\hfill
\begin{equation}
z_k(x,\varrho ):=\left\{
\begin{array}{ll}
y_k(x,\varrho ), & \quad k=0,\ldots ,p-1, \\
y_k(x,\varrho )/\exp (i\varrho \varepsilon _k), & \quad k=p,\ldots ,n-1.
\end{array}
\right.  \label{eq1.10}
\end{equation}
This choice of a fundamental system of solutions is natural due to the fact
that\hfill
\begin{equation}
z_k=\mbox{O}(1)\;,k=0,\ldots ,n-1;\;\;0\le x,\xi \le 1  \label{eq1.11}
\end{equation}
for $\varrho \in S_0, S_1$.\hfill

\subsection{Green's function representation}

Recall that the Green's function admits a representation as a ratio of two
determinants
\begin{equation}
G(x,\xi ,\varrho )= \frac{(-1)^n\Delta (x,\xi ,\varrho )}{n\varrho
^{n-1}\Delta(\varrho )}.  \label{eq1.12}
\end{equation}
Its denominator
\begin{equation}
\Delta (\varrho )=\left| \varrho ^{-j}U_{\jmath }(z_k)\right| _{j,k=0}^{n-1}
\label{eq1.13}
\end{equation}
is usually referred to as \emph{the characteristic determinant}. The
nominator has the form:\hfill
\begin{equation}
\Delta (x,\xi ,\varrho ):=i\cdot \left|
\begin{array}{cc}
z^T & g(x,\xi ,\varrho ) \\
\mathbf{\Delta }(\varrho ) & V(\xi ,\varrho )
\end{array}
\right| ,  \label{eq1.14}
\end{equation}
$z^T$ stands for the row $(z_0(x,\varrho ),\ldots ,z_0(x,\varrho ))$,
\begin{equation}
g(x,\xi ,\varrho ):=g_0(x,\xi ,\varrho )\cdot (n\varrho ^{n-1})/i
\label{eq1.15}
\end{equation}
and
\[
V(\xi ,\varrho )=\left( \varrho ^{-j}U_{jx}(g(x,\xi ,\varrho ))\right)
_{j=0}^{n-1}.
\]
Here the subscript \emph{x} means that the boundary form $U_{\jmath }$ acts
on the kernel $g(x,\xi ,\varrho )$ over the first argument.

Changing a little bit notation from \cite[p.1185]{efs90b} we use an
abbreviation:\hfill
\[
\lbrack [a]]:=a+O(1/\varrho )+ O\Biggl(
\exp  \Bigl(  -\min_{m=0,\ldots ,n-1}|Im(\varrho \varepsilon _m)| \Bigr) %
\Biggr),
\]
where $\varrho $ lies in a fixed sector $S_k,\ \ 0\le k \le n-1$ off some
sufficiently large circle $\{|\varrho |\le R_0\}$.\hfill

Obviously $\mathbf{\Delta }= \left[ {\Delta }_{jk} \right] _{j,k=0}^{n-1}$
is a block-matrix with $r_{\jmath }\times 1$ entries
\begin{equation}
{\Delta }_{jk}:=\varepsilon _k^j\cdot (b_{\jmath }^0\cdot [z_k(0)]+b_{\jmath
}^1\cdot [z_k(1)]).  \label{eq1.16}
\end{equation}
But
\[
\ z_k(1)=[[0]],\,z_k(0)=[[1]],\,\,\,k<p;\quad
z_k(1)=[[1]],\,z_k(0)=[[0]],\,\,\,k\ge p.
\]
Therefore
\begin{equation}
\mathbf{\Delta }(\varrho )= \left[ \left[ \mathbf{\Theta }_p\left( b^0,b^1
\right) \right] \right] .  \label{eq1.17}
\end{equation}
Next, setting\hfill
\begin{eqnarray}
u_k=\left\{
\begin{array}{lcl}
e^{i\varrho \varepsilon _k}\cdot n(\varrho \varepsilon _k)^{n-1}\cdot \tilde
{y_k}(\xi ,\varrho ) & = & e^{i\varrho \varepsilon _k(1-\xi )}\cdot [1], \\
& k<p, &  \\
\tilde {y_k}(\xi ,\varrho )\cdot n(\varrho \varepsilon _k)^{n-1} & = &
e^{i\varrho \varepsilon _k(-\xi )}\cdot [1], \\
& k\ge p &
\end{array}
\right.  \label{eq1.18}
\end{eqnarray}
we come to the following relation\hfill
\begin{eqnarray}
V(\xi ,\varrho ) &=&\quad i\sum_{k<p}[{B}_k^1]\cdot \varepsilon _ku_k(\xi
,\varrho )  \label{eq1.19} \\
&& -i\sum_{k\ge p}[{B}_k^0]\cdot \varepsilon _ku_k(\xi ,\varrho ).  \nonumber
\end{eqnarray}
At last, expanding $\Delta (x,\xi ,\varrho )$ along the $0$-th row we obtain
an important representation:\hfill
\begin{equation}
G(x,\xi ,\varrho )=g_0(x,\xi ,\varrho )+\frac{-2\pi i}{n\varrho ^{n-1}}%
\sum_{t,k=0}^{n-1}a_{tk}(\varrho )\cdot z_k(x,\varrho )\cdot u_t(\xi
,\varrho ),  \label{eq1.20}
\end{equation}
where
\begin{equation}
a_{tk}:=\left\{
\begin{array}{ll}
+\frac{\varepsilon _t}{2\pi }\cdot \left| \Delta \mathop{\longleftarrow}%
\limits_{k} [{B}_t^1]\right| /\Delta , & t<p \\
-\frac{\varepsilon _t}{2\pi }\cdot \left| \Delta \mathop{\longleftarrow}%
\limits_{k} [{B}_t^0],\right| /\Delta , & t\ge p.
\end{array}
\right.  \label{eq1.21}
\end{equation}
Here
\[
\left| \Delta \mathop{\longleftarrow}\limits_{k} d\right|
\]
stands for the determinant $\Delta \left( \varrho \right) $ with the $k$-th
column replaced by a vector $d$. Earlier the coefficients (\ref{eq1.21})
were introduced in \cite{min95a}. Therefore we inserted the factor '$-2\pi i$%
' in (\ref{eq1.21}) in order to preserve notations from that paper as well
as checked a misprint there in the sign.

It would be quite natural to call the matrix ${A}={A}(\varrho)=\left[
a_{tk}\right] _{t,k=0}^{n-1}$ a \textit{characteristic matrix} of the
boundary value problem (\ref{eq1.1})-(\ref{eq1.2}) because it differs from
the analogous object from \cite[p.276]{Nai:ldo} or \cite[p.309]{At:discont}
by another choice of the {fundamental system of solutions}. Namely, in these
books the latter is taken analytic in $\lambda $ with a unit Cauchy data
matrix at the end-point $0$ of the main interval $[0,1]$.

\section{Boundedness of the characteristic matrix\label{sec:bound}}

Further we shall need the following estimate for the resolvent in the lower
half-plane
\begin{equation}
\left\| R_\lambda \right\| \leq \left| Im\lambda \right| ^{-1},Im\lambda <0
\label{eq2.1}
\end{equation}
which stems from the dissipativity of operator $L$. A ray $\left\{ \arg
\lambda =\alpha \right\} ,$ $\pi <\alpha <2\pi $ in $\mathbb{C}_{-}$%
corresponds to the ray $\left\{ \arg \varrho =\alpha /n\right\} $ lying in
the sector $S_1$ in the $\varrho $-plane. Then the estimate (\ref{eq2.1})
turns into
\begin{equation}
\left\| R_\lambda \right\| \leq C\cdot \left| \varrho \right| ^{-n},\quad
\arg \varrho =\alpha /n.  \label{eq2.2}
\end{equation}

\begin{lemma}
\label{lem2.1}The integral operator $g_0$ in $L^2(0,1)$ (see (\ref{eq1.9}) )
admits an estimate

\begin{equation}
\left\| g_0\right\| \leq C\cdot r^{-n},\quad \arg \varrho =\alpha
/n,\;r:=\left| \varrho \right| \geq R_0,  \label{eq2.3}
\end{equation}
where $R_0$ is some fixed positive number.
\end{lemma}

{}

of\textit{\ }Removing square brackets from the asymptotic expressions for
the functions $y_k(x,\varrho )$ and $\tilde y_{\jmath }(\xi ,\varrho )$ we
obtain a kernel $G_0(x,\xi ,\varrho )$ which is naturally extended to $%
\mathbb{R}$. Obviously the latter coincides with the Green's function of the
self-adjoint  operator $D^n$ in $L^2(\mathbb{R}).$ Then
\[
g_0(x,\xi ,\varrho )=G_0(x,\xi ,\varrho )+O\left( \frac 1{\varrho ^n}\right)
.
\]
Clearly operator $G_0$ in $L^2(\mathbb{R})$ with the kernel $G_0(x,\xi
,\varrho )$ obeys estimate (\ref{eq2.3}) which completes the proof.{$\;\;%
\rule{2mm}{2mm}$}

\begin{lemma}
\label{lem2.2}Let $P=P(\varrho )$ be a finite dimensional operator in $%
L^2(0,1)$ with the kernel
\[
P(x,\xi ,\varrho ):=\sum_{t,k=0}^{n-1}a_{tk}(\varrho )\,z_k(x,\varrho
)u_t(\xi ,\varrho ).
\]
Then (\ref{eq2.3}) and representation (\ref{eq1.20}) yield an estimate:
\begin{equation}
\left\| P\right\| _{L^2(0,1)\rightarrow L^2(0,1)}\leq C/r,\mbox{ }r\geq
R_0,\;\arg \varrho =\alpha .  \label{eq2.4}
\end{equation}
In addition, a double-sided estimate holds\hfill
\begin{equation}
\left\| P\right\| _{L^2(0,1)\rightarrow L^2(0,1)}\asymp \sqrt{%
\sum_{t,k=0}^{n-1}\left| a_{tk}(\varrho )\right| ^2\cdot \frac 1{r^2}}
\label{eq2.5}
\end{equation}
uniformly with respect to $\varrho ,r=\left| \varrho \right| \geq R_0,\;\arg
\varrho =\alpha /n$ .
\end{lemma}

{}

of First observe that (\ref{eq2.4}) stems readily from (\ref{eq2.2}),(\ref
{eq2.3}) and (\ref{eq1.20}). Further, let $f\in L^2(0,1).$ Then $%
Pf=\sum_{k=0}^{n-1}d_k\,z_k(x,\varrho )$ with
\[
d_k:=\sum_{t=0}^{n-1}a_{tk}(\varrho )\int_0^1f(\xi )u_t(\xi ,\varrho )\,d\xi
.
\]

Invoking \emph{almost orthogonality }of the f.s.s.\ \ (\ref{eq1.10}), we
arrive at the relation:
\[
\left\| Pf\right\| _{L^2(0,1)}\asymp \sum_{k=0}^{n-1}\left| d_k\right|
^2\,\left\| z_k\right\| _{L^2(0,1)}^2.
\]
Next a direct calculation demonstrates that
\[
\left\| z_k\right\| ^2\asymp \frac 1r,\qquad r\geq R_0,\;\arg \varrho
=\alpha .
\]
Hence it suffices to prove that
\begin{equation}
\sup_{\left\| f\right\| \leq 1}\sum_{k=0}^{n-1}\left| d_k\right| ^2\asymp
\frac 1r\sum_{t,j=0}^{n-1}\left| a_{tj}(\varrho )\right| ^2,\;\forall
k=0,\ldots ,n-1.  \label{eq2.6}
\end{equation}
Fix $\varrho $ and suppose that the sum ($k=0,\ldots,n-1$)
\[
\sum_{t=0}^{n-1}\left| a_{tk}(\varrho )\right| ^2
\]
attains its maximum for $k=k_0(\varrho )$. Then it suffices to check
validity of (\ref{eq2.6}) when $k=k_0.$ The left-hand side of (\ref{eq2.5})
with $k=k_0$ equals
\[
\left\| \sum_{t=0}^{n-1}a_{tk_0}\,u_t(\xi ,\varrho )\right\| _{L^2(0,1)}
\]
which is equivalent to $\sum_{t=0}^{n-1}\left| a_{tk_0}\right|
^2\left\|u_t\right\| ^2$. Here we applied the \emph{almost orthogonality }of
the system (\ref{eq1.18}). But the latter has just an exponential
asymptotics and this is the unique ingredient needed for this property (see
\cite{min93b}). At last, a direct calculation shows that
\[
\left\| u_t\right\| ^2\asymp \frac 1r,\qquad r\geq R_0,\;\arg \varrho
=\alpha/n
\]
which completes the proof. {$\;\;\rule{2mm}{2mm}$}

\begin{corollary}
The \emph{characteristic matrix} is bounded
\begin{equation}
\sum_{t,k=0}^{n-1}\left| a_{tk}(\varrho )\right| ^2=O\left( 1\right) ,\quad
\arg \varrho =\alpha /n,\ \ \left| \varrho \right| \geq R_0.  \label{eq2.7}
\end{equation}
Indeed, one should compare (\ref{eq2.4}) and (\ref{eq2.5}).
\end{corollary}

\section{End of the proof\label{sec:end}}

Let $A_t$ be the $t$-th column of the matrix ${A}=\left( A_0,\ldots
A_{n-1}\right) $. Then $A_t$ satisfies an equation\hfill
\begin{equation}
\mathbf{\Delta }\cdot A_t=\pm \frac{\varepsilon _t}{2\pi }\left[ {B}%
_t^{\#}\right] ,\quad 0\leq t\leq n-1;\quad \#=\left\{
\begin{tabular}{ll}
1, & $t<p$ \\
0, & $t\geq p$.
\end{tabular}
\right.  \label{eq3.1}
\end{equation}
Here and in what follows the sign $^{\prime }+^{\prime }$ corresponds to the
case $t<p$, $^{\prime }-^{\prime }$ --- to the case $t\geq p$.

\begin{lemma}
\label{lem3.1} For every $t\in \left\{ 0,\ldots ,n-1\right\} $ there exists
a vector $\eta _t\in \mathbb{C}^n$ such that
\begin{equation}
\mathbf{\Theta }\left( b^0,b^1\right) \cdot \eta _t=\pm \frac{\varepsilon _t%
}{2\pi }\left[ {B}_t^{\#}\right] .  \label{eq3.2}
\end{equation}
\end{lemma}

{}

ofFix $t\in \left\{ 0,\ldots ,n-1\right\} .$ Using the compactness of the
set of vectors
\[
A_t(\varrho ),\qquad \arg \varrho =\alpha/n,\ \ \left| \varrho \right| \geq
R_0
\]
we obtain an existence of a vector $\eta _t$ such that
\[
\eta _t=\lim_{m\longrightarrow \infty }A_t(\varrho _m)
\]
for some sequence $\varrho _m\rightarrow \infty $. In the meantime the
formula
\begin{equation}
\lim_{m\longrightarrow \infty }\mathbf{\Delta }(\varrho _m)=\mathbf{\Theta }%
_p\left( b^0,b^1\right)  \label{eq3.3}
\end{equation}
stems immediately from (\ref{eq1.17}). Combine the two formulas above and we
are done. {$\;\;\rule{2mm}{2mm}$}

\begin{lemma}
\label{lem3.2} Denote $R(A)$ the image of the matrix $A$. Then
\begin{equation}
R\left( \mathbf{\Theta }\left( b^0,b^1\right) \right) \supset span\left( {B}%
_0^0,\ldots ,{B}_{p-1}^0,{B}_p^1,\ldots ,{B}_{n-1}^1\right) .  \label{eq3.4}
\end{equation}
\begin{equation}
R\left( \mathbf{\Theta }\left( b^0,b^1\right) \right) \supset span\left( {B}%
_0^1,\ldots ,{B}_{p-1}^1,{B}_p^0,\ldots ,{B}_{n-1}^0\right) .  \label{eq3.5}
\end{equation}
\end{lemma}

{}

ofFirstly, apply the matrix $\mathbf{\Theta }_p\left( b^0,b^1\right) $ to
the standard basis in $\mathbb{C}^n$ and get (\ref{eq3.4}). Secondly, (\ref
{eq3.5}) stems from (\ref{eq3.2}) when $t$ runs over $0,\ldots ,n-1.$ {$\;\;%
\rule{2mm}{2mm}$}

\begin{lemma}
\label{lem3.3} The following final relation is valid:
\begin{equation}
R\left( \mathbf{\Theta }_p\left( b^0,b^1\right) \right) =\mathbb{C}^n.
\label{eq3.6}
\end{equation}
\end{lemma}

{}

ofInclusions (\ref{eq3.4})-(\ref{eq3.5}) yield that
\[
R\left( \mathbf{\Theta }_p\left( b^0,b^1\right) \right) \supset span\left( {B%
}_j^0,j=0,\ldots ,n-1;{B}_0^1,j=0,\ldots ,n-1\right) .
\]
Set
\[
{\Psi }=\left( \varepsilon _j^k\right) _{jk=0}^{n-1},\quad \mathbf{Q}=\left(
\mathbf{Q}^0\mathbf{,Q}^1\right) ,
\]
where $\mathbf{Q}^i:=\left( {B}_t^i,t=0,\ldots ,n-1\right) $ is an $n\times n
$ matrix. Then a $(j,k)$-entry of the product $\mathbf{Q}^i\cdot{\Psi }^{*}$
will be a $r_j\times 1$ vector
\[
b_j^i\sum_{t=0}^{n-1}\varepsilon _t^j\cdot \overline{\varepsilon _t^k}%
=b_j^i\cdot n\cdot \delta _{jk},\qquad i=0,1; \qquad j,k=0,\ldots ,n-1,
\]
whence
\begin{equation}
\frac 1n\mathbf{Q}\cdot\Psi^{*}={B}:=\left(
\begin{tabular}{llllll}
b$_0^0$ &  &  & b$_0^1$ &  &  \\
& $\ddots $ &  &  & $\ddots $ &  \\
&  & b$_{n-1}^0$ &  &  & b$_{n-1}^1$%
\end{tabular}
\right) .  \label{eq3.7}
\end{equation}
Clearly, $rank\,{B=}\sum_{j=0}^{n-1}r_j=n$. Therefore $R \left( \mathbf{Q}%
\Psi^{*} \right) = R(B)=\mathbb{C}^n$. {$\;\;\rule{2mm}{2mm}$}

The proof of theorem \ref{thm0.1} is completed. Indeed, we proved that $%
\mathbf{Q}$ is a full range matrix since ${\Psi }$ is invertible. Hence $%
\mathbf{Q}$ is itself invertible. {$\;\;\rule{2mm}{2mm}$}

\section{Concluding remarks}

\subsection{Dissipative boundary conditions}

Consider an operator $L_{ess}$, generated by the simplest expression $D^n$
and the boundary conditions\hfill
\[
b_{\jmath }^0D^jy(0)+b_{\jmath }^1D^jy^{(}1)=0,\;\;j=0,\ldots ,n-1.
\]
It is natural to call $L_{ess}$ \emph{an essential part} of the differential
operator $L$. Obviously, $L\in(R)\Longleftrightarrow L_{ess}\in(R)$.\hfill

In \cite{sal68} it was shown that $L$ is self-adjoint if and only if so is $%
L_{ess}$. Therefore we are able to speak about self-adjointness of the
boundary conditions themselves. Of course, it would be desirable to
establish the same facts in the case of dissipative differential operators
and therefore to derive \textit{regularity of dissipative boundary
conditions themselves } from the theorem \ref{thm0.1}.\hfill

However, the present description of dissipative boundary conditions is
rather abstract \cite{ut73} and requires further clarification.

\subsection{A limit of the characteristic matrix}

As a byproduct of previous considerations we also established a statement
which seems to be of independent interest.

\begin{theorem}
\label{thm4.1} Let $T_\varepsilon ^{\pm }$ be sectors in $\mathbb{C}_{\pm }$
$\left( 0<\varepsilon <\pi \right) $:
\[
T_\varepsilon ^{+}:=\left\{ \lambda \quad |\quad \varepsilon \leq \arg
\lambda \leq \pi -\varepsilon \right\} ,
\]
\[
T_\varepsilon ^{-}:=\left\{ \lambda \quad |\quad \pi +\varepsilon \leq \arg
\lambda \leq 2\pi -\varepsilon \right\}
\]
and $L$ be an $n$-th order differential operator (not necessarily
dissipative), defined by boundary value problem (\ref{eq1.1})-(\ref{eq1.2}).

\begin{enumerate}
\item  Given a sequence $\left\{ \lambda _m\right\} _1^\infty \subset
T_\varepsilon ^{-}$, such that estimate (\ref{eq2.2}) is fulfilled for $%
\lambda =\lambda _m$, we have
\[
\exists \lim_{m\rightarrow \infty }\mathbf{\Delta }(\lambda _m)=:\mathbf{%
\Delta }_\infty =\mathbf{\Theta }_p\left( b^0,b^1\right) ,\quad
\]
\[
\exists \lim_{m\rightarrow \infty }{A}(\lambda _m)=:{A}_\infty ={A}_\infty
(L)=\mathbf{\Theta }_p\left( b^0,b^1\right) ^{-1}\mathbf{\Theta }_p\left(
b^1,b^0\right) \cdot {D}.
\]
Here ${D}:=diag\left( \varepsilon _0,\ldots ,\varepsilon _{p-1},-\varepsilon
_p,\ldots ,-\varepsilon _{n-1}\right) /\left( 2\pi \right) $; all the
matrices don't vanish.

\item  The same relations are valid if such a sequence exists in the sector $%
T_\varepsilon ^{+}$.
\end{enumerate}
\end{theorem}

\begin{remark}
In the odd order case, $n=2q+1,$ take $\varrho \in S_1$. Then theorem \ref
{thm4.1} yields
\begin{equation}
\Theta \left( b^0,b^1\right) \neq 0  \label{eq4.1}
\end{equation}
because $p=q-1$ (see table \ref{table:p}). However, we are not able to
assert the same with respect to another regularity determinant $\Theta
\left( b^1,b^0\right) $ because there is no information concerning the
resolvent's estimate in the upper half-plane $\mathbb{C}_{+}$, equivalently,
$\varrho \in S_0$.
\end{remark}

Note that if (\ref{eq2.2}) were valid than theorem \ref{thm4.1} would imply
\[
0\neq \Theta _p\left( b^0,b^1\right) \equiv \Theta_q\left( b^0,b^1\right),
\]
but the determinant from the right coincides with $\Theta \left(
b^1,b^0\right) $ up to some nonzero multiplicative constant \cite{sal68}.

Add here that (\ref{eq4.1}) implies \emph{half-regularity} of the boundary
conditions (\ref{eq1.2}). This notion has been recently introduced in \cite
{mishu97} and it was shown there that it yields certain information on the
eigenfunctions' and eigenvalues' behaviour.

\subsection{N.Dunford-J.Schwartz' spectrality}

\begin{remark}
\label{rem4.2} In the case, when $L$ is spectral (in N.Dunford-J.Schwartz'
sense \cite{DunSchw:III}) instead of dissipativity, we recently established
existence of appropriate sequences
\[
\left\{ \lambda _m^{\pm }\right\} _1^\infty \subset T_\varepsilon ^{\pm },
\]
such that estimate (\ref{eq2.2}) is fulfilled for $\lambda =\lambda _m^{\pm }
$. Then theorem \ref{thm4.1} implies $L\in (R)$. However, the construction
of these sequences is rather cumbersome and we shall consider it elsewhere.
\end{remark}

Note that the converse was proved in works of G.M.Kesel\'{}man,
V.P.Mikhailov and A.A.Shkalikov (see \cite{kes64,mikh62,shk79a} or
\cite[p.98--99]{Nai:ldo}). More precisely, Birkhoff-regularity yields
spectrality if the boundary conditions are \emph{strongly regular} or
unconditional basicity of eigenfunctions with brackets otherwise.

Thus remark \ref{rem4.2} \textsf{\ finishes classification of spectral two
point {boundary value problem}s } up to a gap between Birkhoff- and strong
regularity.

\subsection{Abstract Birkhoff-regularity}

Perhaps the most striking fact established in the proof of theorem \ref
{thm0.1} is existence of the limit of the \textit{characteristic matrix} ${A}%
(\varrho )$ as $\lambda \rightarrow \infty $ in any sector lying strictly in
the lower half-plane. It is helpful to note here that this limit is an
invariant of Birkhoff-regular problems. Namely, the following statement has
been recently obtained by our student T.Mizrova. The proof is purely
algebraic and will be published elsewhere.

\begin{theorem}[T.Mizrova,1997]
\label{thm:mizrova} Given two differential operators $L_1,$ $L_2\in (R)$,
assume that the limits of their characteristic matrices coincide: $A_\infty
(L_1)=A_\infty (L_2)$. Then their essential parts also coincide: $%
L_{1ess}\equiv L_{2ess}$.
\end{theorem}

Further, we think that there is a close connection between ${A}(\varrho )$
and the B.S.Nagy-C.Foia\c{s}' \textit{characteristic function} of the
dissipative operator $L$ \cite{SNF:harm}. Moreover, we suppose validity of
the following

\begin{conjecture}
Consider an abstract completely continuous operator $G$. Let $\Lambda
:=\left\{ \lambda _j\right\} $ be the set of its eigenvalues and let
\[
K_\lambda :=\left\{ \left| z-\lambda \right| \leq \delta \cdot \left(
1+\left| Im\lambda \right| \right) \right\}
\]
be a circle centered at the point $\lambda $, where $\delta >0$ is some
fixed number. Then existence of an invertible limit of the characteristic
function of $G$ in the sectors $T_\varepsilon ^{\pm }$ off the circles $%
K_{\lambda _j},K_{\overline{\lambda _j}},\lambda _j\in \Lambda $ is a
correct reformulation of the Birkhoff-regularity condition and yields
unconditional basicity of eigenvectors of $G$, perhaps, with brackets.
\end{conjecture}

Recall that in this and more general situations \textit{characteristic
function } was defined by M.S.Liv\v{s}ic, see, for instance, \cite
{Brod:triang}.

\subsection{Quasidifferential expressions}

All results of the paper may be applied as well for a general
quasidifferential expression $l(y)\equiv y^{[n]}$ of the form
\begin{equation}
y^{[0]}=y,\ \ y^{[j]} =Dy^{[j-1]}+\sum_{k=0}^{j-1}p_{j-1,k}(x)y^{[k]},\quad
j=1,\ldots ,n  \label{eq4.2}
\end{equation}
with summable coefficients $p_{j-1,k}$ because the main tool---  property of
\textit{almost orthogonality }--- is valid for such expressions \cite{min93b}%
.



\newpage

\begin{thebibliography}{99}
\bibitem{At:discont}  \textsc{Atkinson F.V.}, \textit{Discrete and
Continuous Boundary Problems} (Mir, Moscow, 1968) (Russian).

\bibitem{Brod:triang}  \textsc{Brodskii M.S.} \textit{Triangular and jordan
representations of linear operators }(Nauka, Moscow, 1969) (Russian).

\bibitem{DunSchw:III}  \textsc{Dunford N.,Schwartz J.T.} \textit{Linear
operators. Spectral operators} (Mir, Moscow, 1974) (Russian).

\bibitem{efs90b}  \textsc{Eberhard W.,Freiling G.,Schneider A.} `{Expansion
theorems for a class of regular indefinite eigenvalue problems}', \textit{J.
Differential Integr. Equat.} 3 (1990), $no^0$ 6, 1181--1200.

\bibitem{fie72}  \textsc{Fiedler H.,} \language=2 `{Zur Re\-gu\-la\-rit\"at
selbst\-ad\-jun\-gier\-ter Rand\-wert\-auf\-gaben}', \language=0 \textit{%
Manuscripta Math.} 7 (1972) 185--196.

\bibitem{kes64}  \textsc{Kesel\'{}man G.M.} 'On unconditional convergence in
eigenfunction expansions for certain differential operators', \textit{Izv.
Vuzov, Matematika} N2(39) (1964) 82--93 (Russian).

\bibitem{hro91a}  \textsc{Khromov A.P.} `Eigenfunction expansion for one
third order boundary value problem', \textit{Matematika i ee prilojeniya} N2
(1991) 17--24 (Saratov univ. press) (Russian).

\bibitem{mikh62}  \textsc{Mikhailov V.P.} 'On Riesz bases in $L_2(0,1)$',
\textit{Dokl. AN SSSR} 144, N5 (1962) 981--984 (Russian).

\bibitem{min77a}  \textsc{Minkin A.M.,} `{Regularity of self-adjoint
boundary conditions}', \textit{Mat. Zametki} 22 (1977), $no^0$ 6, 835--846
(Russian); Engl. transl. \textit{Math. Notes} 22 (1978) 958--965.

\bibitem{min93a}  \textsc{\mbox{\rule{3em}{.4pt}}\thinspace } `{On the class
of regular boundary conditions', \textit{Results Math.} 24 (1993), $n^0$
3/4, 274--279. }

\bibitem{min80}  \textsc{\mbox{\rule{3em}{.4pt}}\thinspace } `Localization
principle for eigen--function series of ordinary differential operators%
\textit{',} \textit{Differenzialniye uravneniya i theoriya functii
(Differential equations and function theory)}. Saratov, Saratov univ. press $%
n^0$3(1980) 68--80 (Russian).\hfill

\bibitem{min93b}  \textsc{\mbox{\rule{3em}{.4pt}}\thinspace } `{A}lmost
orthogonality of {B}irkhoff's {S}olutions', \textit{Results Math.} 24
(1993), $n^0$ 3/4, 280--287.

\bibitem{min95a}  \textsc{\mbox{\rule{3em}{.4pt}}\thinspace } `Odd and Even
cases of Birkhoff-regularity', \textit{Math. Nachr.} 174 (1995) 219--230.

\bibitem{mishu97}  \textsc{\mbox{\rule{3em}{.4pt}}\thinspace and {S}huster L.%
} `Estimates of eigenfunctions for one class of boundary conditions',
\textit{Bull. London Math. Soc.} 29 (1997) 459--469.

\bibitem{SNF:harm}  \textsc{B.Sz.Nagy and C.Foia\c s} \textit{Harmonic
analysis of operators in hilbert space} (Mir, Moscow, 1970) (Russian).

\bibitem{Nai:ldo}  \textsc{M.A. Naimark}, \textit{{L}inear differential
operators} (Nauka, Moscow, 1969) (Russian).

\bibitem{Nik:shift}  \textsc{Nikol\'{}ski\u \i {} N.K.} \textit{Lectures on
the shift operator} (Moscow, Nauka, 1979) (Russian); English transl.,
\textit{Treatise on the shift operator} (Springer, 1986).

\bibitem{sal68}  \textsc{Salaff S.}, `{R}egular {B}oundary {C}onditions for {%
O}rdinary {D}ifferential {O}perators', \textit{Trans. Amer. Math. Soc.} 134
(1968), $no^0$ 2, 355--373.

\bibitem{shk79a}  \textsc{Shkalikov A.A.} 'Basicity of eigenfunctions of an
ordinary differential operator', \textit{Uspekhi mat. nauk} 34, N5(20)
(1979) 235-236 (Russian).

\bibitem{tre:89}  \textsc{Treil\'{} S.R.} '{Geometric methods in spectral
theory of vector-valued functions: some recent results.}' \textit{Operator
Theory: Advances and Applications} (Birkh\"auser Verlag, Basel) 42 (1989)
209--280.

\bibitem{ut73}  \textsc{Utkin V.I.} `Extensions of closed dissipative
operators', \textit{Mat. Zametki} 14, N2 (1973) 223--232 (Russian).
\end{thebibliography}
\end{document}